\documentclass[11pt]{article}
\usepackage{amscd, color}
\usepackage{bbm}
\usepackage{mathrsfs}
\usepackage{bbm}
\usepackage{amsfonts}
\usepackage{amsmath}
\usepackage{stmaryrd}

\usepackage{latexsym,amsfonts,amssymb}

\newtheorem{thm}{Theorem}[section]
\newtheorem{prop}[thm]{Proposition}
\newtheorem{lem}[thm]{Lemma}
\newtheorem{rem}[thm]{Remark}
\newtheorem{coro}[thm]{Corollary}

\begin{document}
\baselineskip=16pt

\title{\bf Twisted vertex operators and  \\
 unitary Lie algebras}

\author{ Fulin Chen$^{3}$, Yun Gao$^{1}$, Naihuan Jing$^{2}$, Shaobin Tan$^{3}$\\
{\small $^{{1}}$ Department of Mathematics and Statistics, York University}\\
 {\small Toronto, Canada M3J 1P3}\\
{\small $^{{2}}$ Department of Mathematics, North Carolina State University }\\
{\small Raleigh, NC, USA 27695}\\
{\small $^{{3}}$ Department of Mathematics, Xiamen University }\\
{\small Xiamen, China 361005}}

\date{}
\maketitle
\begin{abstract}
A representation of the central extension of the 
unitary Lie algebra
coordinated with a skew Laurent polynomial ring
is constructed using vertex operators over an integral $\mathbb Z_2$-lattice.
The irreducible decomposition of the representation is explicitly computed and described.
As a by-product, some fundamental representations of affine Kac-Moody Lie algebra of type $A_n^{(2)}$ are recovered
by the new method.

\end{abstract}

\section{Introduction}\label{section 1} \setcounter{equation}{0}

Affine Kac-Moody Lie algebras, or nontrivial central extensions of loop algebras, are a class of
infinite dimensional Lie algebras fundamentally important in mathematics and theoretical physics. They
are first realized by vertex operators in \cite{LW} and \cite{KKLW} for the principal picture
and later in \cite{FK, S} for the homogeneous picture.
During the last two decades these constructions have been generalized from several directions. The vertex
operator representations of the toroidal Lie algebras have been
given in \cite{MRY, Y, EM} and \cite{T3} (see also \cite{FJW}) in the homogenous setting and
in \cite{B} and \cite{T2}
in the principle setting; vertex representations of quantum affine algebras
have been obtained in \cite{FJ} for untwisted cases and in \cite{J}
for twisted cases; the vertex representations of the
extended affine Lie algebras of type $A$ coordinated by a quantum torus have been provided
in \cite{BS} and \cite{G4} via the principle construction and in \cite{G3} via the homogenous
construction. Finally in \cite{BGT} a unified treatment of vertex representations of affine Lie algebras
using mixed bosons and fermions have been studied, and
the Tits-Kantor-Koecher algebra has also been realized in \cite{T1} using the vertex operator calculus.

On the other hand, before the development of toroidal Lie algebras,
the elementary unitary Lie algebra $\mathfrak{eu}_\nu(\mathcal{R},\bar{ }\ )$ associated with an
involutive associative algebra $\mathcal{R}$ was studied in [AF] as a derived subalgebra
of the unitary Lie algebra $\mathfrak{u}_\nu(\mathcal{R},\bar{ }\ )$. These Lie algebras are generalization of
the usual loop algebras by replacing the commutative coordinate ring with an noncommutative algebra.
As is well-known that the central extensions of Lie algebras are essentially given by
A. Connes' cyclic homology \cite{KL}. For the extended affine Lie algebras the central extensions are also described
by the dihedral homology \cite{G1, ABG}, and the relevant central extensions in the unitary algebras coordinated with
an involutive algebra will be given by a Steinberg Lie algebra structure.

The Steinberg group is the universal central extension of the commutator subgroup of the general linear group.
The Steinberg Lie algebra is
defined as an associative algebra generated by the generators $x_{ij}(\lambda)$ subject to the
Steinberg relations. The extended affine Lie algebras
share many common features with loop algebras and also have distinguished central extensions
given by dihedral homology. When one further generalizes the
algebraic structure by replacing the coordinate ring with the quantum polynomial ring or even an associative algebra,
Steinberg
unitary Lie algebras shall be the key for the universal central extensions. The realization of
Steinberg unitary Lie algebras is thus a natural question in this regard.

In this paper, inspired by the work of Wakimoto in \cite{W}, we will construct a family of twisted
vertex operators associated to an integral $\mathbb{Z}_2$-lattice of rank $\nu$ to realize the
Steinberg relations.
More precisely, for any non-zero complex number $a$,
we define vertex operators $X_{ij}(a,z)$ for $1\le i,j\le \nu$ on a Fock space $M$.
For the purpose of computing the commutator relations, two multi-product decompositions of $\delta$-function and $D\delta$-function are developed.
 For any abelian group $G$ with a character $\sigma: G\to \mathbb{C}^{\times}$,
we can define a corresponding skew Laurent-polynomial ring $\mathcal{R}_\sigma$
and an anti-involution $\bar{ }$ over $\mathcal{R}_\sigma$. We will show that the vertex operators $X_{ij}(a,z)$ give a representation for the Lie algebra
$\widehat{\mathfrak{u}}_\nu(\mathcal{R}_\sigma,\bar{ }\ )$ and also
provide a non-trivial central extension of
the unitary Lie
algebra $\mathfrak{u}_\nu(\mathcal{R}_\sigma,\bar{ }\ )$ associated to the pair $(\mathcal{R}_\sigma, \bar{ }\ )$,
namely, the twisted vertex representation is actually
a module for the central extension $\widehat{\mathfrak{eu}}_\nu(\mathcal{R}_\sigma,\bar{ }\ )$.

As was shown in \cite{W}, it is interesting to analyze the irreducible decomposition of
the twisted module. Based on the observation that the twisted group algebra of the integral $\mathbb{Z}_2$-lattice is isomorphic
to a (finite dimensional) Clifford algebra we are able to use the well-known classical
representation theory of finite dimensional Clifford algebra to achieve this goal.
We remark that this is slightly different from the approach taken by Wakimoto in \cite{W}.
Let $|\sigma|$ be the order of the group $\sigma(G)$, if $|\sigma|=\infty$ or $|\sigma|\in 2\mathbb{N}$, then each
irreducible component of $\widehat{\mathfrak{u}}_\nu(\mathcal{R}_\sigma,\bar{ }\ )$-module $M$ remains irreducible as
$\widehat{\mathfrak{eu}}_\nu(\mathcal{R}_\sigma,\bar{ }\ )$-module.
For the case $|\sigma|\in 2\mathbb{N}+1$, note that in this case all the elements in $\sigma(G)$ are roots of unity.
To this end, we introduce a conjugate anti-involution $\tau$ on
$\widehat{\mathfrak{u}}_\nu(\mathcal{R}_\sigma,\bar{ }\ )$ as all the elements of $\sigma(G)$ lie in the unit circle in the complex plane.
Then, we prove that the $\widehat{\mathfrak{u}}_\nu(\mathcal{R}_\sigma,\bar{ }\ )$-module $M$ is unitary with respect to $\tau$,
which in turn deduce that $M$ is completely reducible as $\widehat{\mathfrak{eu}}_\nu(\mathcal{R}_\sigma,\bar{ }\ )$-module even if $|\sigma|\in 2\mathbb{N}+1$.
As a by-product when taking $G=\{1\}$, we get a completely reducible module $M$ for the affine Kac-Moody Lie algebra of type $A_{\nu-1}^{(2)}$, which recovers the result in \cite{W}.

The paper is organized as follows. In Sect.2, we define a family of vertex
operators for any $a\in \mathbb{C}^{\times}$ and then compute their commutation relations.
In Sect.3, using the derived relations we construct the vertex operator representations for the unitary Lie algebra with non-trivial central extension and then determine the
irreducible components.
Finally, in Sect.4, we construct vertex representations for the elementary unitary
Lie algebra with non-trivial central extension. As an example we also obtain a new realization of the affine Kac-Moody
algebra of type $A_{\nu-1}^{(2)}$.

Throughout this paper, we denote the field of complex numbers,
the group of non-zero complex numbers, the ring of
integers and the set of non-negative integers by $\mathbb{C},\mathbb{C}^{\times},
\mathbb{Z}$ and $\mathbb{N}$, respectively.

\section{Fock Space and Vertex Operators}\label{section 2} \setcounter{equation}{0}
Let $\nu\geq 2$ be a positive integer.
Let $\Gamma=\oplus_{i=1}^{\nu}\mathbb{Z}\epsilon_i,$
 $(\epsilon_i,\epsilon_j)=\delta_{ij}$ for $i,j=1,\cdots, \nu$,
 and $\epsilon_i(n)$ be a linear copy of
 $\epsilon_i$ for $n\in 2\mathbb{Z}+1$ and $i=1,\cdots, \nu$. We define a Lie algebra
 $$
 \mathcal {H}=span_{\mathbb{C}}\{\epsilon_i(n), 1|n\in 2\mathbb{Z}+1\}
 $$
 subject to the following Lie algebra relation
 $$
 [\alpha(m), \beta(n)]={\frac m 2}(\alpha,\beta)\delta_{m+n,0}
 $$
 for $ \alpha, \beta\in\Gamma$, and $m,n\in 2\mathbb{Z}+1$. Let
 ${\bar\Gamma}=\Gamma/2\Gamma$ be the quotient additive group of
 $\Gamma$ factored by the subgroup $2\Gamma$. Let $\mathbb{C}
 [{\bar\Gamma}]=\oplus_{\alpha\in{\bar\Gamma}}\mathbb{C}e^{\alpha}$
 be the twisted group over the finite group $\bar\Gamma$ with
 multiplication defined by
 $e^{\bar\alpha}e^{\bar\beta}=\varepsilon(\alpha,\beta)e^{\overline{\alpha+\beta}}$
 for $\alpha, \beta\in\Gamma$, where the two-cocycle
 $\varepsilon$: $\Gamma\times\Gamma\to \{\pm 1\}$ is given by
\begin{equation}
 \varepsilon(\epsilon_i,\epsilon_j)=
 \begin{cases}1,\ &\text{if}\ i\leq j;\\
 -1,\ &\text{if}\ i>j.
\end{cases}\end{equation}
 and
 $$
 \varepsilon(\sum m_i\epsilon_i, \sum
 n_j\epsilon_j)=\prod_{i,j}(\varepsilon(
 \epsilon_i,\epsilon_j))^{m_in_j}.
 $$
We note that the multiplication in $\mathbb{C}
[{\bar\Gamma}]$ is well-defined as
 $\varepsilon(\alpha_1,\beta_1)=\varepsilon(\alpha_2,\beta_2)$ if
 $\alpha_2-\alpha_1, \beta_2-\beta_1\in 2\Gamma$. Let $z, w$ be
 formal variables and $\alpha\in \Gamma$, set
$$
E_{\pm}(\alpha,z)=\exp\left(-2\sum_{n\in\pm(2\mathbb{N}+1)}{\frac
{\alpha(n)} n}z^{-n}\right),
$$
and define the Fock space
\begin{align}
M=\mathbb{C}[{\bar\Gamma}]\otimes{\mathcal S}({\mathcal H}^-),
\end{align}
where ${\mathcal S}(\mathcal {H}^-)$ is the usual commutative symmetric
algebra over the subalgebra $\mathcal{H}^-$ of $\mathcal{H}$ spanned by
$\epsilon_i(n)$ with $1\le i\le \nu$ and $n\in-(2\mathbb{N}+1).$\\

\begin{lem}\rm  For $\alpha,\beta\in\Gamma$, one has
$$
E_{\pm}(\alpha,z)=E_{\pm}(-\alpha,-z),
$$
$$E_+(\alpha,z)E_-(\beta, w)=E_-(\beta,
w)E_+(\alpha,z)\left({\frac
{1-w/z}{1+w/z}}\right)^{(\alpha,\beta)}.
$$ \\
\end{lem}

\noindent\textbf{Proof.}
The first identity is clear. To prove the second identity,
we set
$$ A=-2\sum_{n\in(2{\mathbb N}+1)}{\frac {\alpha(n)} n}z^{-n},
\;\; B=-2\sum_{n\in-(2{\mathbb N}+1)}{\frac {\beta(n)} n}w^{-n},
$$
and note that $$ [A,B]=4\sum_{m\in(2{\mathbb N}+1)}\sum_{n\in-(2{\mathbb
N}+1)}{\frac {z^{-m}w^{-n}}{mn}}[\alpha(m), \beta(n)]
$$
$$
=4\sum_{m\in(2{\mathbb N}+1)}{\frac
{z^{-m}w^{m}}{-m^2}}(\alpha,\beta){\frac m 2}
$$
$$
=-2(\alpha,\beta)\sum_{m\in(2{\mathbb N}+1)}{\frac 1
m}(w/z)^m=\log\left({\frac
{1-w/z}{1+w/z}}\right)^{(\alpha,\beta)}.
$$
Then by applying the formal rule $e^Ae^B=e^{[A,B]}e^Be^A$ if
$[A,B]$ commutes with $A$ and $B$, we get the required identity.
\hfill$\Box$\\

For $a\in {\mathbb C}^{\times}$, we define the following vertex
operators, which act on the Fock space $M$ as the usual way (see
[FLM]).
$$
X_{ij}(a,z)=
$$
\begin{equation} \begin{cases}
\varepsilon(\epsilon_i,\epsilon_j)e^{\overline{\epsilon_i-\epsilon_j}}E_-(\epsilon_i,
z)E_-(-\epsilon_j, az)E_+(\epsilon_i, z)E_+(-\epsilon_j,
az),\ &i\neq j;\\
 4\epsilon_i(z),\ &i=j, a=1;\\
 \frac{1+a}{1-a}\Big(E_-(\epsilon_i, z)E_-(-\epsilon_i,
az)E_+(\epsilon_i, z)E_+(-\epsilon_i, az)-1\Big),\ &i=j,
a\neq1.
\end{cases}\end{equation}
where $\epsilon_i(z)=\sum_{n\in 2{\mathbb Z}+1}\epsilon_i(n)z^{-n}$,
and $e^{\overline{\epsilon_i-\epsilon_j}}$ is the usual operator
acting on the group algebra ${\mathbb C}[{\bar\Gamma}]$ twisted by the
two-cocycle $\varepsilon$. Set
$$A_{\pm}=-2\sum_{n\in\pm(2{\mathbb N}+1)}{\frac {\epsilon_i(n)}
n}(az)^{-n},\;\; B_{\pm}=-2\sum_{n\in\pm(2{\mathbb N}+1)}{\frac
{\epsilon_i(n)} n}z^{-n},$$
 Then for $a\neq 1$, we have
\begin{align*}
 &{\frac 1 {1-a}}\left( E_{\pm}(\epsilon_i, z)-E_{\pm}(\epsilon_i,
 az)\right)\\
 =&{\frac 1 {1-a}}\left(e^{B_{\pm}}-e^{A_{\pm}}\right)
 ={\frac 1 {1-a}}\left(\sum^{\infty}_{l=1}{\frac 1
 {l!}}B^l_{\pm}-\sum^{\infty}_{l=1}{\frac 1 {l!}}A^l_{\pm}\right)\\
=&{\frac 1 {1-a}}\sum^{\infty}_{l=1}{\frac 1 {l!}}\left[\left(
2\sum_{n\in\pm(2{\mathbb N}+1)}{\frac {\epsilon_i(n)}
n}(az)^{-n}-2\sum_{n\in\pm(2{\mathbb N}+1)}{\frac {\epsilon_i(n)}
n}z^{-n}\right)\sum_{j=0}^{l-1}A_{\pm}^{l-1-j}B_{\pm}^j\right]\\
=&\sum^{\infty}_{l=1}{\frac 1 {l!}}\left[2\sum_{n\in\pm(2{\mathbb
N}+1)}{\frac {\epsilon_i(n)} n}(z)^{-n}{\frac
{a^{-n}-1}{1-a}}\sum_{j=0}^{l-1}A_{\pm}^{l-1-j}B_{\pm}^j\right]\\
\longrightarrow &2\sum_{n\in\pm(2{\mathbb
N}+1)}\epsilon_i(n)z^{-n}\sum_{l=1}^{\infty}{\frac 1
{l!}}lB_{\pm}^{l-1}= 2\sum_{n\in\pm(2{\mathbb
N}+1)}\epsilon_i(n)z^{-n}E_{\pm}(\epsilon_i, z),
\end{align*}
as $a\to 1$.  From this we obtain, for $a\neq 1$, that
\begin{align*}
X_{ii}(a, z)=&{\frac {1+a}{1-a}}E_-(-\epsilon_i, az)\Big(
E_-(\epsilon_i,z)E_+(\epsilon_i,z)-E_-(\epsilon_i,az)E_+(\epsilon_i,az)\Big)E_+(-\epsilon_i,
az)\\
=&({1+a})E_-(-\epsilon_i,
az)\Bigg({\frac{E_-(\epsilon_i,z)-E_-(\epsilon_i,az)}{1-a}}E_+(\epsilon_i,z)+\\
& E_-(\epsilon_i,az){\frac {E_+(\epsilon_i,z)-
E_+(\epsilon_i,az)}{1-a}}\Bigg)E_+(-\epsilon_i,az)\\
\longrightarrow &4\epsilon_i(z)
\end{align*}
as $a \to 1$. Therefore, we have
\begin{align}
\lim_{a\to 1}X_{ii}(a,z)=X_{ii}(1,z).\end{align}

Moreover, from the definition of the vertex operator $X_{ij}(a,
z)$, we can easily check the following result\\

\begin{lem}\rm  For $a\in {\mathbb C}^{\times}$, and $1\leq i,
j\leq \nu$, one has
$$
X_{ij}(a,z)=-X_{ji}(a^{-1}, -az).
$$ \\ \end{lem}

The vertex operator $X_{ij}(a, z)$ can be formally expanded as
follows
$$
X_{ij}(a, z)=\sum_{n\in {\mathbb Z}}x_{ij}(a, n)z^{-n},
$$
where $x_{ij}(a, n)$ are operators acting on the Fock space $M$.
Then the equation in the previous lemma implies that
\begin{align}
x_{ji}(a^{-1}, n)=-(-1)^na^nx_{ij}(a,n),
\end{align}
for $n\in {\mathbb Z}$, $i,j=1,\cdots, \nu$, $a\in {\mathbb C}^{\times}$.

\begin{lem}\rm  For $\alpha\in\Gamma$, set
$\alpha(z)=\sum_{n\in 2{\mathbb Z}+1}\alpha(n)z^{-n}$, then
$$
[\alpha(z), E_{\pm}(\beta, w)]=(\alpha,\beta)E_{\pm}(\beta,
w)\sum_{n\in\pm(2{\mathbb N}+1)}({\frac z w})^n.
$$\end{lem}

\noindent\textbf{Proof.} The result follows from the formal rule $[A, e^B]=[A,B]e^B$
if $[A,B]$ commutes with $B$, and the identity
$$
[\alpha(z), -2\sum_{n\in\pm(2{\mathbb N}+1)}{\frac {\beta(n)}
n}w^{-n}]= (\alpha,\beta)\sum_{n\in\pm(2{\mathbb N}+1)}({\frac z
w})^n.\eqno\Box
$$

The following result is well-known (see [FLM]), and will be used
frequently
later on.\\

\begin{lem}\rm Let $Y(w,z)$  be a formal power series in
$w,z$ with coefficient in a vector space, such that
$\lim_{z\rightarrow w}Y(w,z)$ exists (in the sense of [FLM]). Set
$D_{z}=z\frac{\partial}{\partial z}$, then
\begin{align}
Y(w,z)\delta(a\frac{w}{z})&=Y(w,aw)\delta(a\frac{w}{z}),\\
Y(w,z)(D\delta)(a\frac{w}{z})&=Y(w,aw)(D\delta)(a\frac{w}{z})+(D_{z}Y)(w,z)\delta(a\frac{w}{z}).
\end{align} \\
 \end{lem}

Now we are going to compute the Lie product $[X_{ij}(a,z),
X_{kl}(b,w)]$ for $1\le i,j,k,l\le \nu$, and $a,b\in{\mathbb
C}^{\times}$. For this purpose we need the following several
combinatorial identities
\\

 \begin{lem}\rm Let $B$ and $A_i$ for $i=0,1,\cdots, n$ be non-zero
 distinct complex numbers, then
\begin{align*}
 &\left(\frac{B}{B-A_0}\right)^2\prod^n_{j=1}{\frac B
 {B-A_j}}=
 \sum_{i=1}^n\left(\frac{A_i}{A_i-A_0}\right)^2\left(\prod_{1\le j\le n,j\not=i}{\frac
 {A_i}{A_i-A_j}}\right){\frac B {B-A_i}}\\
 &+\left(\prod_{i=1}^n \frac{A_0}{A_0-A_i}\right)\left[
 \left(\frac{B}{B-A_0}\right)^2+\left(\sum_{l=1}^n \frac{A_l}{A_l-A_0}\right)
 \frac{B}{B-A_0}\right],\\
&\left(\frac{A_0}{B-A_0}\right)^2\prod^n_{j=1}{\frac {A_j}
 {B-A_j}}=
 \sum_{i=1}^n\left(\frac{A_0}{A_i-A_0}\right)^2\left(\prod_{1\le j\le n,j\not=i}{\frac
 {A_j}{A_i-A_j}}\right){\frac {A_i} {B-A_i}}\\
 &+\left(\prod_{i=1}^n \frac{A_i}{A_0-A_i}\right)\left[
 \left(\frac{A_0}{B-A_0}\right)^2+\left(\sum_{l=1}^n \frac{A_0}{A_l-A_0}\right)
 \frac{A_0}{B-A_0}\right].
\end{align*}
 for $n\ge 1$.
 \end{lem}

\noindent\textbf{Proof.} We will only prove the first identity by induction on $n$.
The proof of the second
one is similar to that of the first one, and we omit it.

It is not difficult to check the result for $n=1$. For the case $n+1$, we have
\begin{align*}
&\sum_{i=1}^{n+1}\left(\frac{A_i}{A_i-A_0}\right)^2\left(\prod_{1\le j\le n+1,j\not=i}{\frac
 {A_i}{A_i-A_j}}\right){\frac B {B-A_i}}+\\
 &\left(\prod_{i=1}^{n+1} \frac{A_0}{A_0-A_i}\right)\left[
 \left(\frac{B}{B-A_0}\right)^2+\left(\sum_{l=1}^{n+1} \frac{A_l}{A_l-A_0}\right)
 \frac{B}{B-A_0}\right]\\
 =&\left(\frac{A_{n+1}}{A_{n+1}-A_0}\right)^2
 \left(\prod_{j=1}^n{\frac
 {A_{n+1}}{A_{n+1}-A_j}}\right){\frac B {B-A_{n+1}}}\\
 &+\sum_{i=1}^n
 \left(\frac{A_i}{A_i-A_0}\right)^2\left(\prod_{1\le j\le n,j\not=i}
 {\frac
 {A_i}{A_i-A_j}}\right){\frac
 {A_i}{A_i-A_{n+1}}}{\frac B {B-A_i}}\\
 &+\left(\prod_{i=1}^n \frac{A_0}{A_0-A_i}\right)\left[
 \left(\frac{A_0}{A_0-A_{n+1}}\frac{B}{B-A_0}\right)
 \left( \frac{B}{B-A_0}-\frac{A_{n+1}}{A_0-A_{n+1}}+\sum_{l=1}^{n}
 \frac{A_l}{A_l-A_0}\right)\right]
 \end{align*}
 \begin{align*}
 =&\left(\frac{A_{n+1}}{A_{n+1}-A_0}\right)^2
 \left(\prod_{j=1}^n{\frac{A_{n+1}}{A_{n+1}-A_j}}\right)
 {\frac B {B-A_{n+1}}}\\
 &+\sum_{i=1}^n
 \left(\frac{A_i}{A_i-A_0}\right)^2\left(\prod_{1\le j\le n,j\not=i}
 {\frac
 {A_i}{A_i-A_j}}\right)\left({\frac
 {A_{n+1}}{A_i-A_{n+1}}}+{\frac B {B-A_i}}\right)
 \frac{B}{B-A_{n+1}}\\
 &+\left(\prod_{i=1}^n \frac{A_0}{A_0-A_i}\right)\left[
 \left({\frac B {B-A_0}}\right)^2
 -\left({\frac
 {A_{n+1}}{A_0-A_{n+1}}}\right)^2\right]
 \frac{B}{B-A_{n+1}}\\
 &+\left(\prod_{i=1}^n \frac{A_0}{A_0-A_i}\right)
 \left(\sum_{l=1}^n \frac{A_l}{A_l-A_0}\right)
 \left({\frac{A_{n+1}}{A_0-A_{n+1}}}+{\frac B {B-A_0}}\right)
 \frac{B}{B-A_{n+1}}\\
 =&\left(\frac{A_{n+1}}{A_{n+1}-A_0}\right)^2
 \left(\prod_{j=1}^n{\frac{A_{n+1}}{A_{n+1}-A_j}}\right)
 {\frac B {B-A_{n+1}}}+F_1+F_2,
 \end{align*}
where
\begin{align*}
F_1=&\sum_{i=1}^n
 \left(\frac{A_i}{A_i-A_0}\right)^2\left(\prod_{1\le j\le n,j\not=i}
 {\frac{A_i}{A_i-A_j}}\right)\left({\frac
 {A_{n+1}}{A_i-A_{n+1}}}\right)
 \frac{B}{B-A_{n+1}}\\
 &+\left(\prod_{i=1}^n \frac{A_0}{A_0-A_i}\right)\left[
 \left(\sum_{l=1}^n \frac{A_l}{A_l-A_0}\right){\frac{A_{n+1}}{A_0-A_{n+1}}}
 -\left({\frac
 {A_{n+1}}{A_0-A_{n+1}}}\right)^2\right]
 \frac{B}{B-A_{n+1}},\\
 F_2=&\sum_{i=1}^n
 \left(\frac{A_i}{A_i-A_0}\right)^2\left(\prod_{1\le j\le n,j\not=i}
 {\frac
 {A_i}{A_i-A_j}}\right)\left({\frac B {B-A_i}}\right)
 \frac{B}{B-A_{n+1}}\\
 &+\left(\prod_{i=1}^n \frac{A_0}{A_0-A_i}\right)\left[
 \left({\frac B {B-A_0}}\right)^2+
 \left(\sum_{l=1}^n \frac{A_l}{A_l-A_0}\right)
 {\frac B {B-A_0}}
 \right]
 \frac{B}{B-A_{n+1}}
 \end{align*}

 Now, by using induction on $n$, we find
 \begin{align*}
 F_1=&-\left(\frac{A_{n+1}}{A_{n+1}-A_0}\right)^2
 \left(\prod_{j=1}^n{\frac{A_{n+1}}{A_{n+1}-A_j}}\right)
 {\frac B {B-A_{n+1}}},\\
 F_2=&\left(\frac{B}{B-A_0}\right)^2\prod^{n+1}_{j=1}{\frac B
 {B-A_j}},
 \end{align*}
 as required.
\hfill$\Box$\\

\begin{coro}\rm  Let $A_i\not= A_j$ for $i\not= j$ be
nonzero
 complex numbers, then
 \begin{align*}
 &\left(\frac{1}{1-A_0x}\right)^2\prod^n_{j=1}{\frac 1
 {1-A_jx}}=
 \sum_{i=1}^n\left(\frac{A_i}{A_i-A_0}\right)^2\left(\prod_{1\le j\le n,j\not=i}{\frac
 {A_i}{A_i-A_j}}\right){\frac 1 {1-A_ix}}\\
 &+\left(\prod_{i=1}^n \frac{A_0}{A_0-A_i}\right)\left[
 \frac{A_0x}{(1-A_0x)^2}+\left(1+\sum_{l=1}^n \frac{A_l}{A_l-A_0}\right)
 \frac{1}{1-A_0x}\right],\\
&\left(\frac{A_0x}{1-A_0x}\right)^2\prod^n_{j=1}{\frac {A_jx}
 {1-A_jx}}=
 \sum_{i=1}^n\left(\frac{A_0}{A_i-A_0}\right)^2\left(\prod_{1\le j\le n,j\not=i}{\frac
 {A_j}{A_i-A_j}}\right){\frac {A_ix} {1-A_ix}}\\
 &+\left(\prod_{i=1}^n \frac{A_i}{A_0-A_i}\right)\left[
 \frac{A_0x}{(1-A_0x)^2}+\left(\sum_{l=1}^n \frac{A_0}{A_l-A_0}-1\right)
 \frac{A_0x}{1-A_0x}\right].
 \end{align*} \\
 \end{coro}

  \begin{prop}\rm Let $A_i$ for $0\le i\le n$ be distinct nonzero
 complex numbers, then
\begin{align*}
 &\left(\frac{1}{1-A_0x}\right)^2\prod^n_{i=1}{\frac 1 {1-A_ix}}-(-1)^n
 \left(\frac{A_0^{-1}x^{-1}}{1-A_0^{-1}x^{-1}}\right)^2
 \prod^n_{i=1}{\frac {A_i^{-1}x^{-1}} {1-A_i^{-1}x^{-1}}}\\
 =&\sum_{i=1}^{n}\left[\left(\frac{A_i}{A_i-A_0}\right)^2
 \prod_{1\le j\le n,j\not=i}{\frac{A_i}{A_i-A_j}}\right]\delta(A_ix)\\
 &+\left(\prod_{i=1}^n \frac{A_0}{A_0-A_i}\right)\left[
 (D\delta)(A_0x)+\left(1+\sum_{l=1}^n\frac{A_l}{A_l-A_0}\right)\delta(A_0x)\right]
 \end{align*}
 \end{prop}

\noindent\textbf{Proof.} By applying the second identity in Corollary 2.6, we have
\begin{align*}
&(-1)^{n-1}\left(\frac{A_0^{-1}x^{-1}}{1-A_0^{-1}x^{-1}}\right)^2
 \prod^n_{i=1}{\frac {A_i^{-1}x^{-1}} {1-A_i^{-1}x^{-1}}}\\
=&(-1)^{n-1}\sum_{i=1}^n\left(\frac{A_0^{-1}}{A_i^{-1}-A_0^{-1}}\right)^2\left(\prod_{1\le j\le n,j\not=i}{\frac
 {A_j^{-1}}{A_i^{-1}-A_j^{-1}}}\right){\frac {A_i^{-1}x^{-1}} {1-A_i^{-1}x^{-1}}}\\
 &+(-1)^{n-1}\left(\prod_{i=1}^n \frac{A_i^{-1}}{A_0^{-1}-A_i^{-1}}\right)\left[
 \frac{A_0^{-1}x^{-1}}{(1-A_0^{-1}x^{-1})^2}+\left(\sum_{l=1}^n \frac{A_0^{-1}}{A_l^{-1}-A_0^{-1}}-1\right)
 \frac{A_0^{-1}x^{-1}}{1-A_0^{-1}x^{-1}}\right]\\
=&\sum_{i=1}^n\left(\frac{A_i}{A_i-A_0}\right)^2\left(\prod_{1\le j\le n,j\not=i}{\frac
 {A_i}{A_i-A_j}}\right)\frac{A_i^{-1}x^{-1}}{1-A_i^{-1}x^{-1}}\\
 &-\left(\prod_{i=1}^n \frac{A_0}{A_0-A_i}\right)\left[
 \frac{A_0^{-1}x^{-1}}{(1-A_0^{-1}x^{-1})^2}-\left(1+\sum_{l=1}^n \frac{A_l}{A_l-A_0}\right)
 \frac{A_0^{-1}x^{-1}}{1-A_0^{-1}x^{-1}}\right]
\end{align*}

Therefore, from the first identity in Corollary 2.6, we obtain
\begin{align*}
&\left(\frac{1}{1-A_0x}\right)^2\prod^n_{i=1}{\frac 1 {1-A_ix}}-(-1)^n
 \left(\frac{A_0^{-1}x^{-1}}{1-A_0^{-1}x^{-1}}\right)^2
 \prod^n_{i=1}{\frac {A_i^{-1}x^{-1}} {1-A_i^{-1}x^{-1}}}\\
 =&\sum_{i=1}^n\left(\frac{A_i}{A_i-A_0}\right)^2\left(\prod_{1\le j\le n,j\not=i}{\frac
 {A_i}{A_i-A_j}}\right)\left(\frac{1}{1-A_ix}+\frac{A_i^{-1}x^{-1}}{1-A_i^{-1}x^{-1}}\right)\\
 &+\left(\prod_{i=1}^n \frac{A_0}{A_0-A_i}\right)
 \left(\frac{A_0x}{(1-A_0x)^2}-
 \frac{A_0^{-1}x^{-1}}{(1-A_0^{-1}x^{-1})^2}\right)\\
&+\left(\prod_{i=1}^n \frac{A_0}{A_0-A_i}\right)
\left[\left(1+\sum_{l=1}^n \frac{A_l}{A_l-A_0}\right)
 \left(\frac{1}{1-A_0x}+\frac{A_0^{-1}x^{-1}}{1-A_0i^{-1}x^{-1}}\right)\right]\\
 =&\sum_{i=1}^{n}\left[\left(\frac{A_i}{A_i-A_0}\right)^2
 \prod_{1\le j\le n,j\not=i}{\frac{A_i}{A_i-A_j}}\right]\delta(A_ix)\\
 &+\left(\prod_{i=1}^n \frac{A_0}{A_0-A_i}\right)\left[
 (D\delta)(A_0x)+\left(1+\sum_{l=1}^n\frac{A_l}{A_l-A_0}\right)\delta(A_0x)\right]
 \end{align*}
 \hfill$\Box$

\begin{prop}\rm Let $A_i$ for $1\le i\le n$ be distinct nonzero
 complex numbers, then
 \begin{align*}
 \prod^n_{i=1}{\frac 1 {1-A_ix}}-(-1)^n\prod^n_{i=1}{\frac {A_i^{-1}x^{-1}} {1-A_i^{-1}x^{-1}}}
 =\sum_{i=1}^{n}\left(\prod_{1\le j\le n,j\not=i}{\frac
 {A_i}{A_i-A_j}}\right)\delta(A_ix)
 \end{align*}
\end{prop}

\noindent\textbf{Proof.} By using Proposition 2.7 and the fact that
\begin{align*}
(1-A_0x)^2\left(\frac{1}{1-A_0x}\right)^2=1=(1-A_0x)^2\left(\frac{A_0^{-1}x^{-1}}{1-A_0^{-1}x^{-1}}\right)^2,
\end{align*}
one may get that,
\begin{align*}
&\prod^n_{i=1}{\frac 1 {1-A_ix}}-(-1)^n\prod^n_{i=1}
{\frac {A_i^{-1}x^{-1}} {1-A_i^{-1}x^{-1}}}\\
=&(1-A_0x)^2\sum_{i=1}^{n}\left[\left(\frac{A_i}{A_i-A_0}\right)^2
 \prod_{1\le j\le n,j\not=i}{\frac{A_i}{A_i-A_j}}\right]\delta(A_ix)\\
 &+(1-A_0x)^2\left(\prod_{i=1}^n \frac{A_0}{A_0-A_i}\right)\left[
 (D\delta)(A_0x)+\left(1+\sum_{l=1}^n\frac{A_l}{A_l-A_0}\right)\delta(A_0x)\right]\\
=&\sum_{i=1}^{n}\left(\prod_{1\le j\le n,j\not=i}{\frac
 {A_i}{A_i-A_j}}\right)\delta(A_ix)
  \end{align*}
 where have used Lemma 2.4 for the second equation.  \hfill$\Box$ \\

For convenience, we set
$$
:X_{ij}(a,z)X_{kl}(b,w):=
e^{\overline{\epsilon_i-\epsilon_j+\epsilon_k-\epsilon_l}}
E_-(\epsilon_i, z)E_-(-\epsilon_j, az)E_-(\epsilon_k,w)E_-(-\epsilon_l, bw)
$$
$$
\cdot  E_+(\epsilon_i, z)E_+(-\epsilon_j, az)
E_+(\epsilon_k,w)E_+(-\epsilon_l, bw),
$$
where $1\leq i,j,k,l\leq \nu, a,b\in{\mathbb C}^{\times}$ with the condition that
$i\neq j$ if $a=1$ and $k\neq l$ if $b=1$.
The proof of the following results is straightforward, and is omitted\\

 \begin{lem}\rm For $a,b\in {\mathbb C}^{\times}$ and
$\alpha\in\Gamma$,  one has
\begin{align*}
[X_{ij}(a,z),
e^{\bar\alpha}]&=((-1)^{(\epsilon_i-\epsilon_j,\alpha)}-1)e^{\bar\alpha}X_{ij}(a,z),\\
X_{ij}(a,z)X_{kl}(b,w)&=:X_{ij}(a,z)X_{kl}(b,w):P^{ij}_{kl}
\end{align*} where
$i\neq j$ if $a=1$ and $k\neq l$ if $b=1$, and
$$
P^{ij}_{kl}=\varepsilon(\epsilon_i,\epsilon_j)\varepsilon(\epsilon_k,\epsilon_l)
\varepsilon(\epsilon_i-\epsilon_j,\epsilon_k-\epsilon_l)
\left({\frac {1-w/az} {1+w/az}}\right)^{-\delta_{jk}}
$$
$$
\cdot\left({\frac {1-w/z} {1+w/z}}\right)^{\delta_{ik}}
\left({\frac {1-bw/az}{1+bw/az}}\right)^{\delta_{jl}}
\left({\frac{1-bw/z} {1+bw/z}}\right)^{-\delta_{il}}.
$$ \\ \end{lem}

By symmetry, from the previous lemma and (2.1) we obtain
\\

 \begin{lem}\rm  Let $i\neq j$ if $a=1$ and $k\neq l$ if $b=1$, then
$$
[X_{ij}(a,z),
X_{kl}(b,w)]=:X_{ij}(a,z)X_{kl}(b,w):\Delta^{ij}_{kl},
$$
where
\begin{equation}\begin{aligned}
\Delta^{ij}_{kl}
=&\left(\frac{1+a}{1-a}\right)^{\delta_{ij}}
\left(\frac{1+b}{1-b}\right)^{\delta_{kl}}
\varepsilon(\epsilon_i,\epsilon_j)\varepsilon(\epsilon_k,\epsilon_l)
\varepsilon(\epsilon_i-\epsilon_j,\epsilon_k-\epsilon_l)\\
&\cdot\left({\frac {az+w}{az}}\right)^{\delta_{jk}} \left({\frac
{z-w}{z}}\right)^{\delta_{ik}}
\left({\frac {az-bw}{az}}\right)^{\delta_{jl}} \left({\frac
{z+bw}{z}}\right)^{\delta_{il}}Q^{ij}_{kl},
\end{aligned}\end{equation}
and
\begin{equation*}\begin{aligned}
Q^{ij}_{kl}&=\left({\frac {1}{1+{\frac w z}}}\right)^{\delta_{ik}}
 \left({\frac {1}{1-{\frac w {az}}}}\right)^{\delta_{jk}}
  \left({\frac {1}{1-{\frac {bw} z}}}\right)^{\delta_{il}}
  \left({\frac {1}{1+{\frac {bw} {az}}}}\right)^{\delta_{jl}} \\
  &-\left({\frac {{\frac {z} {w}}}{1+{\frac {z}
{w}}}}\right)^{\delta_{ik}}
 \left({\frac {-{\frac {az} {w}}}{1-{\frac {az} {w}}}}\right)^{\delta_{jk}}
  \left({\frac {-{\frac {z} {bw}}}{1-{\frac {z} {bw}}}}\right)^{\delta_{il}}
  \left({\frac {{\frac {az} {bw}}}{1+{\frac {az}
  {bw}}}}\right)^{\delta_{jl}}.
\end{aligned}\end{equation*}\\ \end{lem}

In the rest of this section, we are going to compute the Lie
product $[X_{ij}(a,z), X_{kl}(b,w)]$ with $1\le i,j,k,l\le \nu$,
and $a,b\in{\mathbb C}^{\times}$, where the vertex operators are
defined by (2.3). For this purpose, we divide the argument into
the following three  cases:

Case one: $i\not= j,$ $k\neq l$;

Case two: $i=j$, $k\neq l$;

Case three: $i=j$, $k=l$.\\

We first consider case one with $i\neq j,$ $k\neq l$.

\begin{prop}\rm  Let $i\not= j,$ $k\not= l$, and
$a, b\in {\mathbb
C}^{\times}$. We have\\

 1. If $i,j,k,l$ are distinct integers, then
$[X_{ij}(a,z), X_{kl}(b,w)]=0$;

 2. If $j=k$, $i\not= l$, then
$$
[X_{ij}(a,z), X_{jl}(b,w)]=2X_{il}(ab, z)\delta({\frac w {az}});
$$

 3. If $j=k$, $i= l$, and $ab\not=1$, then
$$
[X_{ij}(a,z), X_{ji}(b,w)]=2X_{ii}(ab, z)\delta({\frac w
{az}})-2X_{jj}(ab, w)\delta({\frac {bw} {z}})$$ $$+2{\frac
{1+ab}{1-ab}}\left(\delta({\frac w {az}})-\delta({\frac {bw}
{z}})\right);
$$

 4. If $j=k$, $i= l$, and $ab=1$, then
$$
[X_{ij}(a,z), X_{ji}(b,w)]=2(X_{ii}(1, z)-X_{jj}(1,
w))\delta({\frac {bw} {z}})+4(D\delta)({\frac {bw}z}).
$$ \end{prop}

\noindent\textbf{Proof.}  From (2.8), one can easily see the first statement.
For the statement two of the proposition, we have
\begin{align*}
\Delta^{ij}_{jl}&=\varepsilon(\epsilon_i,\epsilon_l)\left({\frac
{az+w}{az}}\right)\left[ {\frac 1 {1-{\frac w{az}}}}+{\frac {\frac
{az}w}{1-{\frac {az}w}}}\right]\\
&=\varepsilon(\epsilon_i,\epsilon_l)\left({\frac
{az+w}{az}}\right)\delta({\frac {w}{az}}).
\end{align*}
Thus, by Lemma 2.10, we get
\begin{align*}
&[X_{ij}(a,z),
X_{jl}(b,w)]=\varepsilon(\epsilon_i,\epsilon_l)\left({\frac
{az+w}{az}}\right)\delta({\frac {w}{az}}):X_{ij}(a,z)X_{jl}(b,w):\\
&=2\varepsilon(\epsilon_i,\epsilon_l)\delta({\frac
{w}{az}})e^{\overline{\epsilon_i-\epsilon_l}} E_-(\epsilon_i,
z)E_-(-\epsilon_l, abz)E_+(\epsilon_i, z)E_+(-\epsilon_l, abz)\\
&= 2X_{il}(ab, z)\delta({\frac w {az}}).
\end{align*}
This gives the statement two. Similarly for statement three, we
have from (2.8) and Proposition 2.8
\begin{align*}
\Delta^{ij}_{ji}&=\left({\frac {az+w}{az}}\right)\left({\frac
{z+bw}{z}}\right)\left[ {\frac 1 {1-{\frac w{az}}}}{\frac 1
{1-{\frac {bw}{z}}}}-{\frac {\frac {az}w}{1-{\frac {az}w}}}{\frac
{\frac {z}{bw}}{1-{\frac {z}{bw}}}}\right]\\
&=\left({\frac {az+w}{az}}\right)\left({\frac
{z+bw}{z}}\right)\left(\frac {1}{1-ab}\delta({\frac w{az}})+
\frac{1}{1-\frac{1}{ab}}\delta({\frac {bw}{z}})\right).
\end{align*}
Thus, from Lemma 2.10, we obtain
\begin{align*}
&[X_{ij}(a,z), X_{ji}(b,w)]\\
=&\left({\frac {az+w}{az}}\right)\left({\frac
{z+bw}{z}}\right)\left(\frac {1}{1-ab}\delta({\frac w{az}})+
\frac{1}{1-\frac{1}{ab}}\delta({\frac {bw}{z}})\right)
:X_{ij}(a,z)X_{ji}(b,w):\\
=&2{\frac {1+ab}{1-ab}}\delta({\frac w{az}}) E_-(\epsilon_i,
z)E_-(-\epsilon_i, abz)E_+(\epsilon_i, z)E_+(-\epsilon_i, abz)\\
&+2{\frac {ab+1}{ab-1}}\delta({\frac {bw}{z}}) E_-(\epsilon_j,
w)E_-(-\epsilon_j, abw)E_+(\epsilon_j, w)E_+(-\epsilon_j, abw)\\
=&2X_{ii}(ab, z)\delta({\frac w {az}})-2X_{jj}(ab, w)\delta({\frac
{bw} {z}}) +2{\frac {1+ab}{1-ab}}\left(\delta({\frac w
{az}})-\delta({\frac {bw} {z}})\right).
\end{align*}

 Finally we prove the statement four of the proposition. In this case we
 have from (2.8)
\begin{align*}
\Delta^{ij}_{ji}&=\left({\frac
{z+bw}{z}}\right)^2\left[\left({\frac 1 {1-{\frac
{bw}{z}}}}\right)^2-\left({\frac {\frac {z}{bw}}{1-{\frac
{z}{bw}}}}\right)^2\right]\\
&=\left({\frac {z+bw}{z}}\right)^2\left({\frac
{bw}{z}}\right)^{-1}(D\delta)({\frac {bw}{z}}).
\end{align*}
Thus, from this and Lemma 2.10, we obtain
\begin{align*}
&[X_{ij}(a,z), X_{ji}(b,w)]\\
=&\left({\frac
{z+bw}{z}}\right)^2\left({\frac
{bw}{z}}\right)^{-1}(D\delta)({\frac
{bw}{z}}):X_{ij}(a,z)X_{ji}(b,w):\\
=&4(D\delta)({\frac {bw}{z}})+4\delta(\frac {bw}{z})\Big(2\epsilon_i(z)-2\epsilon_j(az)\Big)\\
=&2\Big(X_{ii}(1,z)-X_{jj}(1,w)\Big)\delta({\frac
{bw}{z}})+4(D\delta)({\frac {bw}{z}}).
\end{align*}
This completes the proof of the Proposition. \hfill $ \Box$\\

Now we consider case two for $[X_{ij}(a,z), X_{kl}(b,w)]$ with
$i=j$, $k\not= l$. The result will be divided into
three subcases.\\

\begin{prop}\rm  Let $1\le i,k,j\le\nu$, and
$a,b\in{\mathbb C}^{\times}$. We have \\

1. If $i,k,j$ are distinct, then $[X_{ii}(a,z), X_{kj}(b,w)]=0$;

2. If $i\not= j$, and $a=1$, then
$$
[X_{ii}(a,z), X_{ij}(b,w)]=2X_{ij}(b,w)\left(\delta({\frac
{w}{z}})-\delta(-{\frac {w}{z}})\right);
$$

3. If $i\not= j$, and $a\not=1$, then
$$
[X_{ii}(a,z), X_{ij}(b,w)]=2X_{ij}(ab,z)\delta({\frac
{w}{az}})+2X_{ji}(ab^{-1},bz)\delta(-{\frac
{w}{z}}).
$$ \end{prop}

\noindent\textbf{Proof.}
The statement one of the proposition is clear. To prove the
second one, we have
\begin{align*}
&[X_{ii}(a,z), X_{ij}(b,w)]=[4\epsilon_i(z), X_{ij}(b,w)]\\
=&4\Big[\epsilon_i(z),
\varepsilon(\epsilon_i,\epsilon_j)e^{\overline{\epsilon_i-\epsilon_j}}E_-(\epsilon_i,w)E_-(-\epsilon_j,bw)
E_+(\epsilon_i,w)E_+(-\epsilon_j,bw)\Big]\\
=&4\varepsilon(\epsilon_i,\epsilon_j)e^{\overline{\epsilon_i-\epsilon_j}}E_-(\epsilon_i,w)E_-(-\epsilon_j,bw)
E_+(\epsilon_i,w)E_+(-\epsilon_j,bw)\sum_{n\in2{\mathbb Z}+1}({\frac z w})^n\\
=&2X_{ij}(b,w)\left(\delta({\frac {w}{z}})-\delta(-{\frac
{w}{z}})\right),
\end{align*}
where in the second last identity we have used Lemma 2.3. We now
prove the last part of the proposition. By (2.8), and Proposition 2.8,
Lemma 2.10, we have
\begin{align*}
&[X_{ii}(a,z), X_{ij}(b,w)]\\
=&{\frac{1+a}{1-a}}\varepsilon(\epsilon_i,\epsilon_j)\Delta_{ji}^{ii}:X_{ii}(a,z)X_{ij}(b,w):\\
=&{\frac
{1+a}{1-a}}\varepsilon(\epsilon_i,\epsilon_j)\left(1+{\frac
w {az}}\right)\left(1-{\frac w {z}}\right)\left[{\frac 1 {1-{\frac
w {az}}}}{\frac 1 {1+{\frac w {z}}}}+{\frac {\frac z w}{1+{\frac z
w}}}{\frac {\frac {az} w}{1-{\frac {az} w}}}\right]
:X_{ii}(a,z)X_{ij}(b,w):\\
=&{\frac
{1+a}{1-a}}\varepsilon(\epsilon_i,\epsilon_j)\left(1+{\frac
w {az}}\right)\left(1-{\frac w {z}}\right)\left[\frac {1}
{1+a}\delta({\frac w {az}})+\frac {1}{1+\frac {1}
{a}}\delta(-{\frac w z})\right]
:X_{ii}(a,z)X_{ij}(b,w):\\
=&2\varepsilon(\epsilon_i,\epsilon_j)e^{\overline{\epsilon_i-\epsilon_j}}
E_-(\epsilon_i,z)E_-(-\epsilon_j,abz)
E_+(\epsilon_i,z)E_+(-\epsilon_j,abz)\delta({\frac w {az}})\\
&-2\varepsilon(\epsilon_i,\epsilon_j)e^{\overline{\epsilon_i-\epsilon_j}}
E_-(\epsilon_i,-az)E_-(-\epsilon_j,-bz)
E_+(\epsilon_i,-az)E_+(-\epsilon_j,-bz)\delta(-{\frac w {z}})\\
=&2X_{ij}(ab,z)\delta({\frac {w}{az}})-2X_{ij}(ba^{-1},-az)\delta(-{\frac
{w}{z}})\\
=&2X_{ij}(ab,z)\delta({\frac {w}{az}})+2X_{ji}(ab^{-1},bz)\delta(-{\frac
{w}{z}})
\end{align*} \hfill $\Box$

Finally we consider case 3, that is to compute $[X_{ij}(a,z),
X_{kl}(b,w)]$ with $i=j=k= l$,
and $a,b\in{\mathbb C}^{\times}$. From the definition of
$X_{ii}(a,z),X_{ii}(b,w)$, we may assume that $a\neq -1, b\neq -1$.

\begin{prop}\rm  Let $1\le i\le\nu$, and
$a,b\in{\mathbb C}^{\times}$. We have

1. For $a=b=1$ , we have
 $$[X_{ii}(1,z), X_{ii}(1,w)]=4\left((D\delta)({\frac wz})-(D\delta)(-{\frac
 wz})\right);
 $$

2. For $a=1$, $b\not=1$,  then
$$
[X_{ii}(1,z), X_{ii}(b,w)]=2X_{ii}(b,w)\left(\delta({\frac
{w}{z}})-\delta(-{\frac {w}{z}})-\delta({\frac
{bw}{z}})+\delta(-{\frac {bw}{z}})\right)
$$ $$
+2{\frac {1+b}{1-b}}\left(\delta({\frac {w}{z}})-\delta(-{\frac
{w}{z}})-\delta({\frac {bw}{z}})+\delta(-{\frac {bw}{z}})\right);
$$

3. If $a\not= 1$, $b\not=1$, and $a\not=b, ab\not=1$, then
\begin{align*}
&[X_{ii}(a,z), X_{ii}(b,w)]\\
=&2\left(X_{ii}(ab,z)+{\frac
{1+ab}{1-ab}}\right)\delta({\frac
{w}{az}})-2\left(X_{ii}(ab,b^{-1}z)+{\frac
{1+ab}{1-ab}}\right)\delta({\frac {bw}{z}})\\
&+2\left(X_{ii}(ab^{-1},bz)+{\frac
{1+ab^{-1}}{1-ab^{-1}}}\right)\delta(-{\frac
{w}{z}})-2\left(X_{ii}(ab^{-1},z)+{\frac
{1+ab^{-1}}{1-ab^{-1}}}\right)\delta(-{\frac {bw}{az}});
\end{align*}

4. If $a\neq 1, b\neq 1$, and $ab=1, a\neq b$, then
\begin{align*}
&[X_{ii}(a,z), X_{ii}(a^{-1},w)]\\
=&2\left(X_{ii}(1,z)-X_{ii}(1,az)\right)\delta(\frac{w}{az})+
4(D\delta)(\frac{w}{az})\\
&+2\left(X_{ii}(a^2,a^{-1}z)+{\frac
{1+a^2}{1-a^2}}\right)\delta(-{\frac
{w}{z}})-2\left(X_{ii}(a^2,z)+{\frac
{1+a^2}{1-a^2}}\right)\delta(-{\frac {w}{a^2z}});
\end{align*}

5. If $a\neq 1, b\neq 1$, and $a=b, ab\neq 1$, then
\begin{align*}
&[X_{ii}(a,z), X_{ii}(a,w)]\\
=&2\left(X_{ii}(a^2,z)+\frac{1+a^2}{1-a^2}\right)
\delta(\frac{w}{az})-2\left(X_{ii}(a^2,a^{-1}z)
+\frac{1+a^2}{1-a^2}\right)\delta(\frac {aw}{z})\\
&-2\left(X_{ii}(1,az)-X_{ii}(1,z)\right)\delta(-\frac{w}{z})
-4(D\delta)(-\frac{w}{z}).
\end{align*}
\end{prop}

\noindent\textbf{Proof.}
 The first part of the proposition is straightforward. To
prove the second part, we have
$$
[X_{ii}(1,z), X_{ii}(b,w)]
$$
$$=
4\Big[\epsilon_i(z), {\frac
{1+b}{1-b}}E_-(\epsilon_i,w)E_-(-\epsilon_i,bw)E_+(\epsilon_i,w)E_+(-\epsilon_i,bw)\Big].
$$
Then by Lemma 2.3, we have
$$
[X_{ii}(1,z), X_{ii}(b,w)]=2{\frac
{1+b}{1-b}}E_-(\epsilon_i,w)E_-(-\epsilon_i,bw)E_+(\epsilon_i,w)E_+(-\epsilon_i,bw)
$$
$$\cdot \left(\delta({\frac {w}{z}})-\delta(-{\frac {w}{z}})-\delta({\frac
{bw}{z}})+\delta(-{\frac {bw}{z}})\right)
$$
$$
=2\left(X_{ii}(b,w)+{\frac {1+b}{1-b}}\right) \left(\delta({\frac
{w}{z}})-\delta(-{\frac {w}{z}})-\delta({\frac
{bw}{z}})+\delta(-{\frac {bw}{z}})\right).
$$
Next, we prove the third part of the proposition.
Note that in this case, one has $1\ne a^{-1} \ne b\ne a^{-1}b$.
Now, from (2.8) and Proposition 2.8, Lemma 2.10,
we have
\begin{align*}
&[X_{ii}(a,z), X_{ii}(b,w)]\\
=&:X_{ii}(a,z)X_{ii}(b,w):
{\frac {1+a}{1-a}}{\frac {1+b}{1-b}}\big({\frac
{az+w}{az}}\big)\big({\frac {z-w}{z}}\big)\big({\frac
{az-bw}{az}}\big)\big({\frac {z+bw}{z}}\big)\\
&\cdot \Big({\frac {1}{1-{\frac {w}{az}}}}{\frac {1}{1+{\frac
{w}{z}}}}{\frac {1}{1+{\frac {bw}{az}}}}{\frac {1}{1-{\frac
{bw}{z}}}} -{\frac {\frac {az}{w}}{1-{\frac {az}{w}}}}{\frac
{\frac {z}{w}}{1+{\frac {z}{w}}}}{\frac {\frac {az}{bw}}{1+{\frac
{az}{bw}}}}{\frac {\frac {z}{bw}}{1-{\frac {z}{bw}}}}\Big)\\
=&:X_{ii}(a,z)X_{ii}(b,w):
{\frac {1+a}{1-a}}{\frac {1+b}{1-b}}\big({\frac
{az+w}{az}}\big)\big({\frac {z-w}{z}}\big)\big({\frac
{az-bw}{az}}\big)\big({\frac {z+bw}{z}}\big)\\
&\cdot \Big[\frac{1}{(1-ba)(1+a)(1+b)}\delta(\frac{w}{az})
+\frac{1}{(1-a^{-1}b^{-1})(1+a^{-1})(1+b^{-1})}\delta(\frac{bw}{z})\\
&+\frac{1}{(1-ba^{-1})(1+b)(1+a^{-1})}\delta(-\frac{w}{z})
+\frac{1}{(1-b^{-1}a)(1+b^{-1})(1+a)}\delta(-\frac{bw}{az})\Big]\\
=&:X_{ii}(a,z)X_{ii}(b,w):\big[\frac{1+ab}{1-ab}
\delta((\frac{w}{az})-\delta(\frac{bw}{z}))
+\frac{1+ab^{-1}}{1-ab^{-1}}
(\delta(-\frac{w}{z})-\delta(-\frac{bw}{az}))\big]\\
=&2\left(X_{ii}(ab,z)+{\frac
{1+ab}{1-ab}}\right)\delta({\frac
{w}{az}})-2\left(X_{ii}(ab,b^{-1}z)+{\frac
{1+ab}{1-ab}}\right)\delta({\frac {bw}{z}})\\
&+2\left(X_{ii}(ab^{-1},bz)+{\frac
{1+ab^{-1}}{1-ab^{-1}}}\right)\delta(-{\frac
{w}{z}})-2\left(X_{ii}(ab^{-1},z)+{\frac
{1+ab^{-1}}{1-ab^{-1}}}\right)\delta(-{\frac {bw}{az}})
\end{align*}

Finally, we will prove the last part of proposition and the
proof of the fourth part
is similar to that of this case, which we omit.
From (2.8) and Proposition 2.7, Lemma 2.10, one has
\begin{align*}
 &[X_{ii}(a,z), X_{ii}(a,w)]\\
=&:X_{ii}(a,z)X_{ii}(a,w):
\big(\frac {1+a}{1-a}\big)^2\big({\frac
{az+w}{az}}\big)\big({\frac {z-w}{z}}\big)^2\big({\frac {z+aw}{z}}\big)\\
&\cdot \Big[{\frac {1}{1-{\frac {w}{az}}}}\left(\frac {1}{1+\frac
{w}{z}}\right)^2{\frac {1}{1-{\frac
{aw}{z}}}} -{\frac {\frac {az}{w}}{1-{\frac {az}{w}}}}\left({\frac
{\frac {z}{w}}{1+{\frac {z}{w}}}}\right)^2{\frac {\frac {z}{aw}}{1-{\frac {z}{aw}}}}\Big]\\
=&:X_{ii}(a,z)X_{ii}(a,w):
\big(\frac {1+a}{1-a}\big)^2\big({\frac
{az+w}{az}}\big)\big({\frac {z-w}{z}}\big)^2\big({\frac {z+aw}{z}}\big)\\
&\Big[\frac{1}{(1+a)^2(1-a^2)}\delta(\frac{w}{az})+\frac{1}{(1+a^{-1})^2(1-a^{-2})}
\delta(\frac{aw}{z})\\
&+\frac{1}{(1+a)(1+a^{-1})}[(D\delta)(-\frac{w}{z})+2\delta(-\frac{w}{z})]\Big]
\end{align*}
\begin{align*}
=&:X_{ii}(a,z)X_{ii}(a,w):\frac{1+a^2}{1-a^2}
\delta((\frac{w}{az})-\delta(\frac{aw}{z}))
-2\left(\epsilon_i(az)-\epsilon_i(z)\right)\delta(-\frac{w}{z})
-4(D\delta)(-\frac{w}{z})\\
&+:X_{ii}(a,z)X_{ii}(a,w):
\frac {a}{(1-a)^2}\Big[2\big({\frac
{az+w}{az}}\big)\big({\frac {z-w}{z}}\big)^2\big({\frac {z+aw}{z}}\big)\\
&-\frac{w}{az}\big({\frac {z-w}{z}}\big)^2\big({\frac {z+aw}{z}}\big)
-\frac{aw}{z}\big({\frac
{az+w}{az}}\big)\big({\frac {z-w}{z}}\big)^2\\
&-2\big(-\frac{w}{z}\big)
\big({\frac
{az+w}{az}}\big)\big({\frac {z-w}{z}}\big)\big({\frac {z+aw}{z}}\big)\Big]
\delta(-\frac{w}{z})\\
=&2\left(X_{ii}(a^2,z)+\frac{1+a^2}{1-a^2}\right)
\delta(\frac{w}{az})-2\left(X_{ii}(a^2,a^{-1}z)
+\frac{1+a^2}{1-a^2}\right)\delta(\frac {aw}{z})\\
&-2\left(X_{ii}(1,az)-X_{ii}(1,z)\right)\delta(-\frac{w}{z})
-4(D\delta)(-\frac{w}{z}).
\end{align*} \hfill $\Box$

\begin{rem}\rm It is easy to see that
$$
{\frac
 {\delta(z)-a^2\delta(az)}{1-a}}\to(D\delta)(z)+2\delta(z)
 $$
 as $a\to 1$. This gives the following identities
$$
{\frac
 {\delta(w/az)-a^2b^2\delta(bw/z)}{1-ab}}\to(D\delta)(w/az)+2\delta(w/az)
 $$
 as $ab\to 1$, and
 $$
{\frac {a^2\delta(-w/z)-b^2\delta(-bw/az)}{a-b}}\to
a(D\delta)(-w/z)+2a\delta(-w/z),
$$
as $b\to a$.
\end{rem}

In summary,
from
Proposition 2.11, 2.12 and 2.13, we can get the following commutator relations

\begin{thm}\rm For $1\le i,j,k,l\le \nu$ and $ a,b\in \mathbb{C}^{\times}$, one has
\begin{align*}
&[X_{ij}(a,z),X_{kl}(b,w)]\\
=&2\delta_{jk}\left(X_{il}(ab,z)+\delta_{il}(1-\delta_{ab,1})
\frac{1+ab}{1-ab}\right)\delta(\frac{w}{az})\\
-&2\delta_{il}\left(X_{kj}(ab,b^{-1}z)+\delta_{jk}
(1-\delta_{ab,1})\frac{1+ab}{1-ab}\right)\delta(\frac{bw}{z})\\
+&2\delta_{ik}\left(X_{lj}(ab^{-1},bz)+\delta_{jl}(1-\delta_{a,b})
\frac{1+ab^{-1}}{1-ab^{-1}}\right)\delta(-\frac{w}{z})\\
-&2\delta_{jl}\left(X_{ik}(ab^{-1},z)+\delta_{ik}
(1-\delta_{a,b})\frac{1+ab^{-1}}{1-ab^{-1}}\right)\delta(-\frac{bw}{az})\\
+&4\delta_{il}\delta_{jk}\delta_{ab,1}(D\delta)(\frac{w}{az})
-4\delta_{ik}\delta_{jl}\delta_{a,b}(D\delta)(-\frac{w}{z}).
\end{align*}
\end{thm}

\section{Representations of unitary Lie algebras}\label{section 3} \setcounter{equation}{0}

In this section we begin by recalling the unitary Lie algebra $\mathfrak{u}_\nu(\mathcal{R}_\sigma,\bar{}\ )$
associated with a skew Laurent polynomial ring $\mathcal{R}_\sigma$ and an anti-involution $\ \bar{ }
\ $, which was first introduced in [AF]. Using the commutator relation among the vertex operators developed in Theorem
2.15, we find that $M$ turn to a representation for a non-trivial central extension of the unitary Lie algebra
$\mathfrak{u}_\nu(\mathcal{R}_\sigma,\bar{}\ )$. Moreover, we determine the irreducible decomposition of
$M$ explicitly.

Let $G$ be an abelian group with a character $\sigma$, that is,
$\sigma:G\rightarrow \mathbb{C}^{\times}$ is a group homomorphism.
Extend $\sigma$, still call $\sigma$,
to be an automorphism of the group algebra $\mathcal{R}$ associated to $G$
determined by $\sigma(e^\alpha)=\sigma(\alpha)e^\alpha, \alpha\in G$.
Then we can form the skew Laurent polynomial ring
$\mathcal{R}_\sigma=\mathcal{R}[t^{\pm1},\sigma]$
with basis $t^me^\alpha, m\in \mathbb{Z},\alpha\in G $ and multiplication
\begin{align*}
e^\alpha t^m=(\sigma(\alpha))^mt^me^\alpha.
\end{align*}
For simplicity, we denote $\sigma(\alpha)$ by $\tilde{\alpha}$ in the following.

 Let $\bar{ }$ be an anti-involution of $\mathcal{R}_\sigma$ defined by
\begin{align*} \bar{t}=-t,\ \overline{e^\alpha}=e^{-\alpha},
\end{align*}
then $\overline{t^me^\alpha}=(-\tilde{\alpha})^{-m}t^me^{-\alpha}$
for $\alpha\in G,m\in \mathbb{Z}$.

Define an operator on the $\nu\times \nu$ matrix algebra $M_\nu(\mathcal{R}_\sigma)$,
\begin{align*}
 ^*:\ M_\nu(\mathcal{R}_\sigma)&\rightarrow M_\nu(\mathcal{R}_\sigma),\\
           X&\mapsto \overline{X}^t
\end{align*}
where $X\in M_\nu(\mathcal{R}_\sigma)$ and $X^t$ is the transpose of
the matrix $X$. Then $\theta(X)=-X^*$ induces an involution on the
Lie algebra $gl_\nu(\mathcal{R}_\sigma)=M_\mu(\mathcal{R}_\sigma)^-$.
The fixed point subalgebra of $\theta$ is called a unitary Lie algebra (see [AF],[G2]):
\begin{align*}
\mathfrak{u}_\nu(\mathcal{R}_\sigma,\bar{ }\ )
=\{X\in M_\nu(\mathcal{R}_\sigma)| X^*=-X\}.
\end{align*}

For convenience, we set
\begin{align*}
e_{ij}(m,\alpha)=E_{ij}t^me^\alpha-E_{ji}\overline{t^me^\alpha},
\end{align*}
where $1\leq i,j\leq \nu, \alpha\in G, m\in \mathbb{Z}$. Then
$e_{ij}(m,\alpha)=-(-\tilde{\alpha})^{-m}e_{ji}(m,-\alpha)$ and
$\mathfrak{u}_\nu(\mathcal{R}_\sigma,\bar{ }\ )$ is spanned by
$e_{ij}(m,\alpha)$ for all $1\leq i,j\leq \nu, m\in \mathbb{Z}$ and
$\alpha\in G$.

We define a 1-dimensional central extension of
$\widehat{gl}_\nu(\mathcal{R}_\sigma)=gl_\nu(\mathcal{R}_\sigma)\oplus
\mathbb{C}c$ of $gl_\nu(\mathcal{R}_\sigma)$ with Lie bracket (See [G3,4])
\begin{equation}\begin{aligned}
&[E_{ij}t^me^\alpha,E_{kl}t^ne^\beta]=
\delta_{jk}\tilde{\alpha}^nE_{il}t^{m+n}e^{\alpha+\beta}\\
&-\delta_{il}\tilde{\beta}^mE_{kj}t^{m+n}e^{\alpha+\beta}
+m\delta_{il}\delta_{jk}\delta_{m+n,0}
\delta_{\tilde{\alpha}\tilde{\beta},1}\tilde{\alpha}^nc
\end{aligned}\end{equation}
where $1\leq i,j,k,l\leq \nu, m,n\in \mathbb{Z}$ and $\alpha,\beta\in G$.

Now we can from a central extension of $\mathfrak{u}_\nu(\mathcal{R}_\sigma,\bar{ }\ )$
\begin{align*}
\widehat{\mathfrak{u}}_\nu(\mathcal{R}_\sigma,\bar{ }\ )=\mathfrak{u}_\nu(\mathcal{R}_\sigma,\bar{ }\ )\oplus \mathbb{C}c
\end{align*}
 with Lie bracket as (3.1). Set
 \begin{align*}
 e_{ij}(\alpha,z)=\sum_{n\in \mathbb{Z}}e_{ij}(n,\alpha)z^{-n},
 \end{align*}
 for $1\leq i,j\leq \nu, \alpha\in G$. Then we have

 \begin{prop}\rm In $\widehat{\mathfrak{u}}_\nu(\mathcal{R}_\sigma,\bar{ }\ )$,
\begin{align*}
&[e_{ij}(\alpha,z),e_{kl}(\beta,w)]\\
=&\delta_{jk}e_{il}(\alpha+\beta,z)\delta(\frac{w}{\tilde{\alpha}z})
-\delta_{il}e_{kj}(\alpha+\beta,\tilde{\beta}^{-1}z)\delta(\frac{\tilde{\beta}w}{z})\\
+&\delta_{ik}e_{lj}(\alpha-\beta,\tilde{\beta}z)\delta(-\frac{w}{z})
-\delta_{jl}e_{ik}(\alpha-\beta,z)\delta(-\frac{\tilde{\beta}w}{\tilde{\alpha}z})\\
+&2\delta_{il}\delta_{jk}\delta_{\tilde{\alpha}\tilde{\beta},1}c(D\delta)(\frac{w}{\tilde{\alpha}z})
-2\delta_{ik}\delta_{jl}\delta_{\tilde{\alpha},\tilde{\beta}}c(D\delta)(-\frac{w}{z}).
\end{align*}
where $1\leq i,j,k,l\leq \nu, \alpha,\beta\in G$.
\end{prop}

\noindent\textbf{Proof.} We have
\begin{align*}
&[\sum_{m\in \mathbb{Z}}(E_{ij}t^me^\alpha-E_{ji}\overline{t^me^\alpha})z^{-m},
\sum_{n\in \mathbb{Z}}(E_{kl}t^ne^\beta-E_{lk}\overline{t^ne^\beta})w^{-n}]\\
=&\sum_{m,n\in \mathbb{Z}}\delta_{jk}(E_{il}t^me^\alpha t^ne^\beta
-E_{li}\overline{t^ne^\beta}\ \overline{t^me^\alpha})z^{-m}w^{-n}\\
-&\sum_{m,n\in \mathbb{Z}}\delta_{il}(E_{kj}t^ne^\beta t^me^\alpha
-E_{jk}\overline{t^me^\alpha} \overline{t^ne^\beta})z^{-m}w^{-n}\\
+&\sum_{m,n\in \mathbb{Z}}\delta_{ik}(E_{jl}\overline{t^ne^\beta}t^me^\alpha
-E_{lj}\overline{t^me^\alpha}t^ne^\beta)z^{-m}w^{-n}\\
-&\sum_{m,n\in\mathbb{Z}}\delta_{jl}(E_{ik}t^me^\alpha\overline{t^ne^\beta}
-E_{ki}t^ne^\beta\overline{t^me^\alpha})z^{-m}w^{-n}\\
+&\sum_{m\in \mathbb{Z}}2m\delta_{jk}\delta_{il}\delta_{\tilde{\alpha}\tilde{\beta},1}
(\tilde{\alpha})^{-m}cz^{-m}w^m
-\sum_{m\in \mathbb{Z}}2m\delta_{ik}\delta_{jl}\delta_{\tilde{\alpha},\tilde{\beta}}c
(-z)^{-m}w^m\\
=&\sum_{m,n\in \mathbb{Z}}\delta_{jk}e_{il}(m+n,\alpha+\beta)z^{-m-n}(\tilde{\alpha}z)^nw^{-n}\\
-&\sum_{m,n\in \mathbb{Z}}\delta_{il}e_{kj}(m+n,\alpha+\beta)
(\tilde{\beta}^{-1}z)^{-m-n}z^n(\tilde{\beta}w)^{-n}\\
+&\sum_{m,n\in \mathbb{Z}}\delta_{ik}e_{lj}(m+n,\alpha-\beta)
(\tilde{\beta}z)^{-m-n}z^n(-w)^{-n}\\
-&\sum_{m,n\in \mathbb{Z}}\delta_{jl}e_{ik}(m+n,\alpha-\beta)
z^{-m-n}(\tilde{\alpha}z)^n(-\tilde{\beta}w)^{-n}\\
+&2\delta_{il}\delta_{jk}\delta_{\tilde{\alpha}\tilde{\beta},1}c(D\delta)(\frac{w}{\tilde{\alpha}z})
-2\delta_{ik}\delta_{jl}\delta_{\tilde{\alpha},\tilde{\beta}}c(D\delta)(-\frac{w}{z}),
\end{align*}
as required.
\hfill $\Box$ \\

Comparing Theorem 2.15 and Proposition 3.1, one can see that the following result holds true.

\begin{thm}\rm Let $G$ be an abelian group with a character $\sigma$, then $M$ is a
vertex operator representation for the unitary Lie algebra $\widehat{\mathfrak{u}}_\nu(\mathcal{R}_\sigma,\bar{ }\ )$
with the action given by
\begin{align*}
2e_{ij}(m,\alpha)
\mapsto x_{ij}(\tilde{\alpha},m)+\delta_{ij}(1-\delta_{\tilde{\alpha},1})\frac{1+\tilde{\alpha}}
{1-\tilde{\alpha}},
\ c\mapsto \frac{1}{2},
\end{align*}
 where $1\leq i,j\leq \nu$ and $\alpha\in G$.
\end{thm}

Next, we are going to determine the irreducible decomposition of
$\widehat{\mathfrak{u}}_\nu(\mathcal{R}_\sigma,\bar{ }\ )$-module $M$.
Let $\mathbb{C}[\bar{Q}]$ be the
subalgebra of $\mathbb{C}[\bar{\Gamma}]$ generated by $e^{\overline{\epsilon_j-\epsilon_{j+1}}}, 1\leq j\leq \nu-1$
which has dimension $2^{\nu-1}$.
Using the standard techniques developed in the vertex representation of affine Lie algebras,
it is easy to see that if $U$ is an irreducible
 $\mathbb{C}[\bar{Q}]$-submodule of $\mathbb{C}[\bar{\Gamma}]$,
 then $U\otimes \mathcal{S}(\mathcal{H}^-)$  remains irreducible as a $\widehat{\mathfrak{u}}_\nu(\mathcal{R}_\sigma,\bar{ }\ )$-module.
 Therefore, we need to determine the irreducible components of $\mathbb{C}[\bar{\Gamma}]$ for detail.
As a by-product, we also give the irreducible decomposition of $\widehat{\mathfrak{u}}_\nu(\mathcal{R}_\sigma,\bar{ }\ )$-module $M$.

We observe that $e^{\bar{\epsilon}_j}e^{\bar{\epsilon}_k}+
e^{\bar{\epsilon}_k}e^{\bar{\epsilon}_j}=2\delta_{jk}, 1\leq j,k\leq \nu$.
So, consider the Clifford algebra $Cl(\nu)$ with generators $\psi_j, 1\leq j\leq \nu$ and
relation
\begin{align*}
\psi_j\psi_k+\psi_k\psi_j=2\delta_{jk},\ 1\leq j,k\leq \nu.
\end{align*}
Since $dim (\mathbb{C}[\bar{\Gamma}])=dim(Cl(\nu))=2^{\nu}$,
$\mathbb{C}[\bar{\Gamma}]$ is isomorphic to
$Cl(\nu)$ as associative algebra.
Then $\mathbb{C}[\bar{Q}]$ is isomorphic to the subspace of
$Cl(\nu)$ spanned by products of an even number of elements $\psi_j, j=1,\cdots,\nu$.

By the classical representation theory for Clifford algebra, we see that
$\mathbb{C}[\bar{Q}]$ is semisimple. Moreover, if $\nu=2d+1$, $\mathbb{C}[\bar{Q}]$
has a unique simple module $S$, which has dimension $2^d$. If $\nu=2d$,
$\mathbb{C}[\bar{Q}]$ has exactly two simple module $S^\pm$,
which has dimension $2^{d-1}$.

\begin{rem}\rm It is well-known that the subspace $Cl_2(\nu)$ of $Cl(\nu)$ spanned
by products of 2
different element $\psi_i,i=1,\cdots,\nu$ is closed under the Lie bracket and
is isomorphic to the simple orthogonal Lie algebra $\mathfrak{o}(\nu)$.
Fix bases in the root systems of types $D_d$ and $B_d$ as follows:
\begin{align*}
 \Pi(D_d)&=\{h_1-h_2,\cdots,h_{d-1}-h_d,h_{d-1}-h_{d} \},\\
 \Pi(B_d)&=\{h_1-h_2,\cdots,h_{d-1}-h_d,h_d\}.
 \end{align*}
Then $S$ is the spinor representation of $\mathfrak{o}(2d+1)$ with
highest weight $\frac{1}{2}(h_1+\cdots+h_d)$. And $S^\pm$ are
the semispinor representation of $\mathfrak{o}(2d)$ with highest
weight $\frac{1}{2}(h_1+\cdots+h_{d-1}\pm h_d)$.
\end{rem}

For $\rho=\pm 1$,  let $w_{2j-1}(\rho)=1+\sqrt{-1}\rho e^{\bar{\epsilon}_{2j-1}}$ and
$w_{2j}(\rho)=1+\rho e^{\bar{\epsilon}_{2j}}$. For $\rho_1,\cdots,\rho_\nu\in \{\pm1\}$,
we put
\begin{align}
w(\rho_1,\cdots,\rho_\nu):=\prod_{j=1}^\nu w_j(\rho_j).
\end{align}
Then the collection of these elements $\{w(\rho_1,\cdots,\rho_\nu): \rho_1,\cdots,\rho_\nu=\pm 1\}$
form a basis of $\mathbb{C}[\bar{\Gamma}]$.
The construction of such basis elements are motivated by the work \cite{W}.
It is easy to check that (or see Lemma 4.1 \cite{W})

\begin{lem}\rm For $\rho_1,\cdots,\rho_\nu\in \{\pm1\}$ and $1\leq j\leq \nu$, we have
\begin{equation*}e^{\bar{\epsilon}_j}w(\rho_1,\cdots,\rho_\nu)=
\begin{cases}\sqrt{-1}\rho_jw(-\rho_1,\cdots,-\rho_j,\rho_{j+1},\cdots,\rho_\nu),\
&\text{if}\ j\ \text{is odd};\\
\rho_j w(-\rho_1,\cdots,-\rho_{j-1},\rho_j,\cdots,\rho_\nu),\
&\text{if}\ j\ \text{is even}.
\end{cases}\end{equation*}
\end{lem}

Then for $1\leq j\leq [\frac{\nu}{2}]$,
we have
\begin{equation}\begin{split}
 e^{\overline{\epsilon_{2j-1}-\epsilon_{2j}}}w(\rho_1,\cdots,\rho_\nu)&=
 -\sqrt{-1}\rho_{2j-1}\rho_{2j}w(\rho_1,\cdots,\rho-\nu),\\
 e^{\overline{\epsilon_{2j}-\epsilon_{2j+1}}}w(\rho_1,\cdots,\rho_\nu)&=
-\sqrt{-1}\rho_{2j}\rho_{2j+1}w(\rho_1,\cdots,\rho_{2j-1},-\rho_{2j},-\rho_{2j+1},
\rho_{2j+2},\cdots,\rho_\nu).
\end{split}\end{equation}
 From (3.3), we see that
the signs $\rho_1,\rho_{2j}\rho_{2j+1}, 1\leq j\leq [\frac{\nu}{2}]$
and $\rho_\nu$ if $\nu$ is even are
invariant under the action of $\mathbb{C}[\bar{Q}]$.

In the case $\nu=2d+1$, for any tuple $\gamma=(\gamma_1,\cdots,\gamma_{d+1})\in \mathbb{Z}_2^{d+1},
\mathbb{Z}_2=\{\pm1\}$, set
\begin{align}
V_\nu(\gamma)=span_\mathbb{C}\{w(\rho_1,\cdots,\rho_\nu)|
 \rho_1=\gamma_1,\rho_{2j}\rho_{2j+1}=\gamma_{j+1},1\leq j\leq d\}.
\end{align} Then $V_\nu(\gamma)$ is a $\mathbb{C}[\bar{Q}]$-module with dimension $2^d$, hence is simple.
 So, we have shown that
 \begin{align*}
 \mathbb{C}[\bar{\Gamma}]=\oplus_{\gamma\in \mathbb{Z}_2^{d+1}}V_\nu(\gamma)
 \end{align*} where each $V_\nu(\gamma)$ is an irreducible $\mathbb{C}[\bar{Q}]$-module.

In the case $\nu=2d$, for any tuple $\gamma=(\gamma_1,\cdots,\gamma_{d+1})\in \mathbb{Z}_2^{d+1}$,
set
\begin{align}V_\nu(\gamma)=span_\mathbb{C}\{w(\rho_1,\cdots,\rho_\nu)|
\rho_1=\gamma_1, \rho_{2j}\rho_{2j+1}=\gamma_{j+1},1\leq j\leq d-1,\rho_{2d}=\gamma_{d+1}\}.
\end{align}
Then $V_\nu(\gamma)$ is an irreducible $\mathbb{C}[\bar{Q}]$-module
with
\begin{align*}
\mathbb{C}[\bar{\Gamma}]=\oplus_{\gamma\in\mathbb{Z}_2^{d+1}}V_\nu(\gamma).
 \end{align*}

In summary, we have obtained
\begin{thm}\rm The $\widehat{\mathfrak{u}}_\nu(\mathcal{R}_\sigma,\bar{ }\ )$-module $M$ is completely reducible.
Moreover,
$$ M=\oplus_{\gamma\in \mathbb{Z}_2^{[\frac{\nu}{2}]+1}}V_\nu(\gamma)\otimes
\mathcal{S}(\mathcal{H}^-),$$
where each submodule $V_\nu(\gamma)\otimes
\mathcal{S}(\mathcal{H}^-)$ is irreducible as $\widehat{\mathfrak{u}}_\nu(\mathcal{R}_\sigma,\bar{ }\ )$-module
and $V_\nu(\gamma)$ are defined in
(3.4) and (3.5). \\
\end{thm}

\section{Application to elementary unitary Lie algebras}\label{section 4} \setcounter{equation}{0}

Let $G$ be an abelian group with a character $\sigma$ and $|\sigma|$ be the order of $\sigma(G)$.
Define the set $\Lambda_\sigma=\{(m,\alpha)\in (\mathbb{Z},G)| m\in |\sigma|\mathbb{Z},
\alpha\in \ker(\sigma)\}$. Here, we set $|\sigma|\mathbb{Z}=0$ if $|\sigma|=\infty$.
Then we have the following basic result about $\mathcal{R}_\sigma$ whose proof is straightforward.

\begin{lem}\rm The derived algebra $\mathcal{R}'_\sigma$ of $\mathcal{R}_\sigma$
has a basis consisting of monomials $t^me^\alpha$ for $(m,\alpha)\not\in \Lambda_\sigma$
 and the center $Z_\sigma$ of $\mathcal{R}_\sigma$ has basis consisting of monomials $t^me^\alpha$
 for $(m,\alpha)\in \Lambda_\sigma$. Therefore, we have $\mathcal{R}_\sigma=\mathcal{R}'_\sigma\oplus Z_\sigma$.
\end{lem}

The elementary unitary Lie algebra $\mathfrak{eu}_\nu(\mathcal{R}_\sigma,\bar{ }\ )$ is the derived algebra of
$\mathfrak{u}_\nu(\mathcal{R}_\sigma,\bar{ }\ )$ which is generated by $e_{ij}(m,\alpha)$ for $1\le i\ne j\le \nu,
m\in \mathbb{Z}$ and $\alpha\in G$.
By applying Lemma 4.1, one can see that
\begin{align*}
\mathfrak{u}_\nu(\mathcal{R}_\sigma,\bar{ }\ )=\mathfrak{eu}_\nu(\mathcal{R}_\sigma,\bar{ }\ )
\oplus \sum_{(n,\beta)\in \Lambda_\sigma}\mathbb{C}(\sum_{i=1}^{\nu}e_{ii}(n,\beta)).
\end{align*}
In particular, we have that
$e_{ii}(m,\alpha)-e_{jj}(m,\alpha),e_{ii}(n,\beta)\in \mathfrak{eu}_\nu(\mathcal{R}_\sigma,\bar{ }\ )$
for $1\le i\ne j\le \nu,m,n\in \mathbb{Z},\alpha,\beta\in G$ and $(n,\beta)\not\in \Lambda_\sigma$.
Let $\widehat{\mathfrak{eu}}_\nu(\mathcal{R}_\sigma,\bar{ }\ )$ be the subalgebra of
$\widehat{\mathfrak{u}}_\nu(\mathcal{R}_\sigma,\bar{ }\ )$ generated by $e_{ij}(m,\alpha)$
for $1\le i\ne j\le \nu, m\in \mathbb{Z}, \alpha\in G$.
Then, $\widehat{\mathfrak{eu}}_\nu(\mathcal{R}_\sigma,\bar{ }\ )$
is a one-dimensional non-trivial central extension of $\mathfrak{eu}_\nu(\mathcal{R}_\sigma,\bar{ }\ )$.

\begin{rem}\rm It was known that
both $\mathfrak{eu}_\nu(\mathcal{R}_\sigma,\bar{ }\ )$ and
$\widehat{\mathfrak{eu}}_\nu(\mathcal{R}_\sigma,\bar{ }\ )$ are $BC_d$-graded with grading subalgebra
$\mathfrak{o}(\nu)$ (See [ABG] for details), where $d=[\frac{\nu}{2}]$.\\
\end{rem}

Next, we turn to consider the unitarizability of $M$.
For our purpose, we may assume that all the elements in $\sigma(G)$ lie in the unity circle in the complex plane,
that is, $|a|=1,\ \forall a\in \sigma(G)$.
Define a conjugate linear map
$\tau:\widehat{\mathfrak{u}}_\nu(\mathcal{R}_\sigma,\bar{ }\ )\rightarrow
\widehat{\mathfrak{u}}_\nu(\mathcal{R}_\sigma,\bar{ }\ )$ by
letting
\begin{align*}
\tau(e_{ij}(m,\alpha))=(-1)^{(m+1)}e_{ij}(-m,\alpha),\ \tau(c)=c,
\end{align*}
where $1\le i,j\le \nu, m\in \mathbb{Z}$ and $\alpha\in G$.
Due to $|\tilde{\alpha}|=1,
\alpha\in G$, one can get that $\tau$ is a conjugate anti-involution of
 $\widehat{\mathfrak{u}}_\nu(\mathcal{R}_\sigma,\bar{ }\ )$.

Define a Hermitian form $<\ |\ >$ on $\mathcal{H}^-$ by
\begin{align*}
<\alpha(m)|\beta(n)>=\frac{m}{2}(\alpha,\beta)\delta_{m,n}, \ \forall\ \alpha(m),\beta(n)
\in \mathcal{H}^-,
\end{align*}
and extend this form to $S(\mathcal{H}^-)$ as usual:
\begin{align*}
<x_1\cdots x_m|y_1\cdots y_n>=\delta_{m,n} \sum_{\sigma\in S_m}\prod_{i=1}^m<x_i|y_{\sigma(i)}>.
\end{align*}
Next, we define $<\ |\ >$ on $\mathbb{C}[\bar{\Gamma}]$ by
\begin{align*} <e^{\bar{\alpha}}|e^{\bar{\beta}}>=\delta_{\bar{\alpha},\bar{\beta}}, \
\forall\ \alpha,\beta\in \Gamma.
\end{align*}

This gives a positive definite Hermitian form on $M$ (See also \cite{FK}(2.12)-(2.14)) such that
\begin{align*} <\epsilon_i(m).w|v>&=<w|\epsilon_i(-m).v>,\\
               <e^{\bar{\epsilon}_i}.w|v>&=<w|e^{\bar{\epsilon}_i}.v>,
\end{align*}
for $1\le i\le \nu, m\in 2\mathbb{Z}+1, w,v\in M$.

\begin{prop}\rm Suppose that all the elements in $\sigma(G)$ lie in the unity circle in the complex plane,
then the $\widehat{\mathfrak{u}}_\nu(\mathcal{R}_\sigma,\bar{ }\ )$-module
$M$ is unitary with respect to $\tau$.
\end{prop}

\noindent\textbf{Proof.} It is sufficient to show that
\begin{align*}
<X_{ij}(a,z).w|v>=<w|-X_{ij}(a,-z^{-1}).v>,
\end{align*}
for $1\le i,j\le \nu, a\in \sigma(G), w,v\in M$.

By applying the fact that $\bar{a}=a^{-1}$, one may get that
\begin{align*}
<\frac{\epsilon_i(m)}{m}(az)^{-m}.w|v>=
<w|\frac{-\epsilon_i(-m)}{-m}(az^{-1})^{m}.v>,
\end{align*}
and that
\begin{align*}
<E_\pm(\epsilon_i,az).w|v>=<w|E_\mp(-\epsilon_i,az^{-1}).v>
=<w|E_\mp(\epsilon_i,-az^{-1}).v>
\end{align*}
for $1\le i\le \nu, a\in \sigma(G),w,v\in M$.

Hence, for $i\ne j$, we find that
\begin{align*}
&<X_{ij}(a,z).w|v>\\
=&<w|E_-(-\epsilon_j,-az^{-1})E_-(\epsilon_i,-z^{-1})
E_+(-\epsilon_j,-az^{-1})E_+(\epsilon_i,-az^{-1})
e^{\bar{\epsilon_j}}e^{\bar{\epsilon_i}}.v>\\
=&<w|-X_{ij}(a,-z^{-1}).v>.
\end{align*}

The case for $i=j, a=1$ is clear.
For the case $i=j, a\ne 1$, one only needs to note that $\overline{\left(\frac{1+a}{1-a}\right)}=
-\frac{1+a}{1-a}$.
\hfill $\Box$ \\

By taking restriction, we know that each component $V_\nu(\gamma)\otimes S(\mathcal{H}^-),
\gamma\in \mathbb{Z}_2^{[\frac{\nu}{2}]+1}$ of $M$ is an
$\widehat{\mathfrak{eu}}_\nu(\mathcal{R}_\sigma,\bar{ }\ )$-module. Moreover,

\begin{thm}\rm Let $V_\nu(\gamma)\otimes S(\mathcal{H}^-)$ be the
$\widehat{\mathfrak{u}}_\nu(\mathcal{R}_\sigma,\bar{ }\ )$-module described in Theorem 3.5,

(1) If $|\sigma|=\infty$ or $|\sigma|\in 2\mathbb{N}$, then $V_\nu(\gamma)\otimes S(\mathcal{H}^-)$ is irreducible as
$\widehat{\mathfrak{eu}}_\nu(\mathcal{R}_\sigma,\bar{ }\ )$-module.

(2) If $|\sigma|=N\in 2\mathbb{N}+1$, let $\mathcal{H}_0^-$ be the subspace of $\mathcal{H}^-$ spanned by
$\epsilon_i(m)-\epsilon_j(m),\epsilon_i(n)$ for $1\le i\ne j\le \nu, m,n\in -(2\mathbb{N}+1), n\not\in N\mathbb{Z}$,
then $V_\nu(\gamma)\otimes S(\mathcal{H}_0^-)$ is an irreducible
$\widehat{\mathfrak{eu}}_\nu(\mathcal{R}_\sigma,\bar{ }\ )$-submodule of $V_\nu(\gamma)\otimes S(\mathcal{H}^-)$.
Moreover, $V_\nu(\gamma)\otimes S(\mathcal{H}^-)$ is completely reducible
as $\widehat{\mathfrak{eu}}_\nu(\mathcal{R}_\sigma,\bar{ }\ )$-module.
\end{thm}

\noindent\textbf{Proof.} First, we write
\begin{align*}
L=\widehat{\mathfrak{eu}}_\nu(\mathcal{R}_\sigma,\bar{ }\ ),\ V=V_\nu(\gamma)\otimes S(\mathcal{H}^-),\
V_0=V_\nu(\gamma)\otimes S(\mathcal{H}_0^-),
\end{align*}
for simplicity.

If $|\sigma|=\infty$ or $|\sigma|\in 2\mathbb{N}$, then $e_{ii}(m,0)\in
L, 1\le i\le \nu, m\in 2\mathbb{Z}+1$ and hence $V$
remains irreducible as $L$-module.

Next, if $|\sigma|=N\in 2\mathbb{N}+1$, we observe that
\begin{align*}
E_-(\epsilon_i,z)E_-(-\epsilon_j,az)\in End(V_0)[[z]],
\end{align*}
for $1\le i,j\le \nu$ and $a\in \sigma(G)$.
From the definition of vertex operators and the fact that
$e_{ii}(m,0)\not\in L,
m\in 2\mathbb{Z}+1\cap N\mathbb{Z}$, we have obtained that the subspace
$V_0$ is invariant under the action of
$L$.
As the usual reason, one can easily see that
$V_0$ is irreducible as
$L$-module.

Finally, we will prove that
$V_\nu(\gamma)\otimes S(\mathcal{H}^-)$ is completely reducible
as $\widehat{\mathfrak{eu}}_\nu(\mathcal{R}_\sigma,\bar{ }\ )$-module if $|\sigma|\in 2\mathbb{N}+1$.
Let $d$ be the degree derivation on $\mathcal{R}_\sigma$ such that
$d.t^me^\alpha=mt^me^\alpha$ for $m\in \mathbb{Z},\alpha\in G$.
Let $\tilde{L}$ be the semi-direct product of the Lie algebra $L$ and the derivation $d$.
Extend $\tau$ to a conjugate anti-involution of $\tilde{L}$ by setting $\tau(d)=d$.
Let $d_0$ be the degree operator on $S(\mathcal{H}^-)$ determined by
$d_0.\epsilon_i(m)=m\epsilon_i(m)$ for $1\le i\le \nu, m\in -(2\mathbb{N}+1)$.
Then, by letting the action of $d$ as $1\otimes d_0$, $V$ can be extended to an
$\tilde{L}$-module and is unitary with respect to $\tau$.
One has a natural $\mathbb{Z}$-graded structure on $\tilde{L}$ and $V$:
\begin{align*}
\tilde{L}=\oplus_{n\in \mathbb{Z}}\tilde{L}_n,\ V=\oplus_{m\in -\mathbb{N}}V(m),
\end{align*}
with respect to the action of $d$. For any submodule $V'$ of $V$,
write $V'(m)=V(m)\cap V'$ for $m\in -\mathbb{N}$ so that $V'=\oplus_m V'(m)$.
Note that each graded subspace $V(m)$ is finite dimensional.
Therefore, each $\tilde{L}$-submodule of $V$ has a orthogonal compliment being unitary
(in this case, each $a\in \sigma(G)$ is a root of unity so that $V$ is unitary by Proposition 4.2).

Let $W$ be the sum of all irreducible $\tilde{L}$-submodules (which contain $V_0$) and
$W'$ be its orthogonal compliment.
Let $W'=\oplus_{m\le n'}W'(m)$ with $W'(n')\ne 0$ for some $n'\in -\mathbb{N}$.
Let $U$ be the $\tilde{L}_0$-module generated by $W'(n')$, then $U\subseteq V(n')$ and is finite
dimensional. Notice that $\tau(\tilde{L}_0)=\tilde{L}_0$,
then $U$ is a finite dimensional unitary
$\tilde{L}_0$-module and hence completely reducible.
Let $U_0$ be an irreducible component of $\tilde{L}_0$-module $U$
and consider the $\tilde{L}$-module $U'$ generated by $U_0$.
Now, $U'$ is an irreducible $\tilde{L}$-submodule of $W'$ as $U'(n')=U_0$ can't split.
But $W$ is supposed to contain all the irreducible $\tilde{L}$-submodules of $V$,
a contradiction.
Therefore, $V=W$ and $V$ is completely reducible as $\tilde{L}$-module (also as $L$-module).
 \hfill $\Box$ \\

In the rest of this paper, we will pay our attention to the trivial case that $G=\{1\}$.
In this case, the elementary unitary Lie algebra $\widehat{\mathfrak{eu}}_\nu(\mathcal{R}_\sigma,\bar{ }\ )$
is isomorphic to the twisted affine
Kac-Moody algebra of type $A_{\nu-1}^{(2)}$.
Let $\Lambda_j,0\le j\le d, d=[\frac{\nu}{2}]$ be the fundamental integral weight of $A_{\nu-1}^{(2)}$.
Following from Theorem 4.4 and the techniques developed by Wakimoto (See \cite{W} Theorem 4.2) , we can find that

\begin{coro}\rm  $M$ is a completely reducible representation
for the affine Kac-Moody algebra of type $A_{\nu-1}^{(2)}$ with the action given in Theorem 3.2.
 If $\nu=2d$, then $V_{2d}(\gamma)\otimes S(\mathcal{H}_0^-),
\gamma\in \mathbb{Z}_2^{d+1}$ are the irreducible $\widehat{\mathfrak{eu}}_\nu(\mathcal{R}_\sigma,\bar{ }\ )$-modules with
highest weight $\Lambda_{d-1}$ if $\gamma_1\cdots\gamma_{d+1}=-1$ or $\Lambda_d$ if $\gamma_1\cdots\gamma_{d+1}=1$.
If $\nu=2d+1$, then $V_{2d+1}(\gamma)\otimes S(\mathcal{H}_0^-),
\gamma\in \mathbb{Z}_2^{d+1}$ are the irreducible $\widehat{\mathfrak{eu}}_\nu(\mathcal{R}_\sigma,\bar{ }\ )$-modules
with highest weight
$\Lambda_d$.\\
\end{coro}

\begin{rem}\rm
The vertex operator representations
$\mathbb{C}[\bar{Q}]\otimes S(\mathcal{H}_0^-)$ for affine Kac-Moody algebra of type $A_{\nu-1}^{(2)}$
have been studied in \cite{W}, where the author used a different two cocycle $\varepsilon(\ ,\ ):Q\to Q$ from our's.
Explicitly, let $\Pi=\{\alpha_i=\epsilon_i-\epsilon_{i+1},1\le i\le \nu-1\}$ be the simple root system of type $A_{\nu-1}$,
then the two-cocycle $\varepsilon(\ ,\ )$ is determined by $\varepsilon(\alpha_i,\alpha_j)(1\le i,j\le \nu-1)$ and its
bi-multiplicative property. The Dynkin diagram of $\Pi$ with orientation corresponds to $\varepsilon$ as
follows:
\begin{equation*}
\varepsilon(\alpha_i,\alpha_j)=\begin{cases}
1,\ &\text{if}\ \begin{CD}
\substack{\alpha_i \\ \bigcirc} @>>> \substack{\alpha_j \\ \bigcirc}
\end{CD}\ \text{or}\ \alpha_i\ \text{is not connected with}\ \alpha_j,\\
-1,\ &\text{if}\ i=j \ \text{or}\ \begin{CD}
\substack{\alpha_i \\ \bigcirc} @<<< \substack{\alpha_j \\ \bigcirc}
\end{CD}.
\end{cases}\end{equation*}
The orientation of Dynkin diagram used in \cite{W} was
\begin{align*}
\nu=2d:\qquad
&\substack{\alpha_1 \\ \bigcirc} \longrightarrow \substack{\alpha_2 \\ \bigcirc}
\longleftarrow \substack{\alpha_3 \\ \bigcirc}
\longrightarrow \substack{\alpha_4 \\ \bigcirc}\longleftarrow \substack{\alpha_5\\ \bigcirc} \cdots \cdots
\substack{\alpha_{2d-3} \\ \bigcirc} \longrightarrow \substack{\alpha_{2d-2} \\ \bigcirc}
\longleftarrow \substack{\alpha_{2d-1} \\ \bigcirc}\\
\nu=2d+1:\qquad &\substack{\alpha_1 \\ \bigcirc} \longrightarrow \substack{\alpha_2 \\ \bigcirc}
\longleftarrow \substack{\alpha_3 \\ \bigcirc}
\longrightarrow \substack{\alpha_4 \\ \bigcirc}\longleftarrow \substack{\alpha_5 \\ \bigcirc} \cdots \cdots
\substack{\alpha_{2d-2} \\ \bigcirc} \longleftarrow \substack{\alpha_{2d-1} \\ \bigcirc}
\longrightarrow \substack{\alpha_{2d} \\ \bigcirc}
\end{align*}
while the orientation used in this paper as follows
\begin{align*}
\substack{\alpha_1 \\ \bigcirc} \longrightarrow \substack{\alpha_2 \\ \bigcirc}
\longrightarrow \substack{\alpha_3 \\ \bigcirc}
\longrightarrow \substack{\alpha_4 \\ \bigcirc}\longrightarrow \substack{\alpha_5 \\ \bigcirc} \cdots \cdots
\substack{\alpha_{\nu-3} \\ \bigcirc} \longrightarrow \substack{\alpha_{\nu-2} \\ \bigcirc}
\longrightarrow \substack{\alpha_{\nu-1} \\ \bigcirc}
\end{align*}
\end{rem}

\begin{rem}\rm Let $G=\mathbb{Z}$ and consider
 the subalgebra $\mathcal{G}$ of $\widehat{\mathfrak{eu}}_\nu(\mathcal{R}_\sigma,\bar{ }\ )$
generated by the elements $e_{ij}(0,m), 1\le i\ne j\le \nu, m\in \mathbb{Z}$. Then,
$\mathcal{G}$ is a fixed point subalgebra of
affine Kac-Moody algebra $A_{\nu-1}^{(1)}$.
But, $\mathcal{G}$ is not isomorphic to the (untwisted or twisted) affine Kac-Moody algebra.
The representation theory of $\mathcal{G}$ are known little so far.
It is interesting to see that  $\tau(x)=-x$ for all $x\in \mathcal{G}$ if $|\sigma(1)|=1$.
\end{rem}

\end{document}